%Authors: Y.Benyamini

%Title: The Uniform Classification of Banach Spaces

%Filename: benyaminiunifrm.tex
%TeX: AMSTeX
%Length: 69407
%Received Date: 6/7/94
%SubjectClass: 46B20
%Abstract: This is a survey of results on the classification of Banach
%spaces as
%metric spaces. It is based on a series of lectures I gave at the
%Functional Analysis Seminar in 1984-1985, and it appeared in the
%1984-1985
%issue of the Longhorn Notes. I keep receiving requests for copies,
%because
%some of the material here does not appear elsewhere and because the
%Longhorn Notes are not so easy to get.  Having it posted on the
%Bulletin
%thus seems reasonable despite the fact that it is not updated, and I
%thank
%the Editors of the Longhorn Notes for their permission to do so.
%

%Citation: Texas functional analysis seminar 1984--1985 (Austin, Tex.),
%15--38, Longhorn Notes, Univ. Texas Press, Austin, TX, 1985.

%32   space        33 ! exclam. pt.   34 " double quote  35 # sharp
%36 $ dollar       37 % percent       38 & ampersand     39 ' prime
%40 ( left paren.  41 ) rt. paren.    42 * asterisk      43 + plus
%44 , comma        45 - minus         46 . period        47 / division
%58 : colon        59 ; semi-colon    60 < less than     61 = equal
%62 > greater than 63 ? question mark 64 @ at
%91 [ left bracket 92 \ backslash     93 ] right bracket 94 ^ caret
% 95 _ underline    96 ` left single quote
%123 { left brace  124 | vertical bar 125 } right brace  126 ~ tilda

% -- Paper for Longhorn Notes 1984-5 --  Author: Prof. Benyamini
%% 5/25/94: TeX copy e-mailed to Benyamini    
%%	 mar29aa%technion.bitnet@cunyvm.cuny.edu
\input amstex
\documentstyle{amsppt}
	\overfullrule=0pt
	\def\IR{{\Bbb R}}
	
\topmatter
\title The Uniform Classification of Banach Spaces\endtitle 
\author Yoav Benyamini\endauthor
\thanks Partially supported by NSF Grant DMS 8403669.\endthanks 
\affil The Technion\\ 
Israel Institute of Technology\\
Haifa 32000 Israel\\
and\\
The University of Texas \\
Austin, Texas 78712\endaffil
\abstract 
This is a survey of results on the classification of Banach spaces as 
metric spaces. It is based on a series of lectures I gave at the 
Functional Analysis Seminar in 1984-1985, and it appeared in the 1984-1985 
issue of the Longhorn Notes. I keep receiving requests for copies, because 
some of the material here does not appear elsewhere and because the 
Longhorn Notes are not so easy to get.  Having it posted on the Bulletin 
thus seems reasonable despite the fact that it is not updated, and I thank 
the Editors of the Longhorn Notes for their permission to do so. 
\endabstract 
\endtopmatter

\document 

\head \S0.  Introduction.\endhead 

Banach spaces are topological (metric) spaces with an additional
structure --- the vector space structure.  In the linear theory
we study all these structures simultaneously, and we deal with
linear continuous maps.  Two spaces are identified if they are
linearly homeomorphic.  One could, however, consider Banach spaces
as a special class of topological (or metric) spaces and study them as 
such.  In the topological classification two spaces are identified
if they are homeomorphic.  In the metric classification we
identify uniformly-homeomorphic spaces.

While the linear theory is very rich --- there are many different
types of spaces, with complicated subspace structure, the
topological theory is, in some sense, trivial.  A remarkable
theorem of Kadec says that any two separable Banach spaces are
homeomorphic.  (See \cite{BP} for a thorough study of Banach
spaces as topological spaces.)  Kadec's theorem was 
extended by Torunczyk \cite{T}, who proved that two Banach spaces
are homeomorphic iff they have the same density character.

The uniform theory lies between these two extremes.  It is 
rich enough so that to say that two Banach spaces are uniformly
homeomorphic already says something about similarities in
their linear structure.  Yet it does not, in general, imply
linear isomorphism.

In this series of lectures I presented some of the ideas and
results in the theory of uniform classification.  There is no
attempt at a comprehensive survey.  I chose one topic --- the
infinite dimensional classification problem, and presented, with
complete proofs, what I believe to be the main ideas and results
on this subject.  Even within this narrower subject I did not
try to be exhaustive.  For example, I didn't even mention when a 
result applies to more general spaces (linear metric spaces,
Polish groups) --- and many of the results do.

There are no new results here, but some of the proofs have
been simplified a little (notably \S 6).

Let me mention briefly two topics on which I didn't speak and
on which some significant work has been done.

The first is the problem of uniform embedding, which is very
similar in spirit to the subject of these lectures.  The
following theorem of Aharoni \cite{Ah} (see also \cite{As} 
for a simpler proof) dramatically shows the difference between
the linear and uniform theories.  While $c_0$ is, linearly,
a very ``small'' space, we have

\proclaim{Theorem 0.1}  Every separable metric space Lipschitz
embeds into $c_0$.
\endproclaim 

The second result, due to Aharoni, Maurey and Mityagin \cite{AMM}
gives a complete solution to the question which spaces embed 
uniformly into a Hilbert space.

\proclaim{Theorem 0.2}  A linear metric space embeds uniformly
into a Hilbert space iff it is linearly isomorphic to a
subspace of $L^0$ --- the space of measurable functions with
the topology of convergence in measure.\endproclaim

In particular $L^p$ embeds uniformly into a Hilbert space
iff $p\leq 2.$

For earlier results on the impossibility of uniform embedding
of some spaces (including $c_0$) into a Hilbert space, see
the very short and elegant papers of Enflo \cite{E1, E2}
as well as his \cite{E3}.  The first two papers already
belong to the second subject I would like to mention --- the 
local theory, or the theory of finite metric spaces.  Rather
than quote some of the results here, let me just say that
there is a close similarity, in spirit and in methods, between
this theory and the theory of finite dimensional Banach spaces,
and refer the reader to few of the recent papers \cite{L-J}, 
\cite{Bo}, \cite{B-M-W} and their references.

We now briefly describe the content of these lectures.

The first four sections deal with Lipschitz classification 
and Lipschitz maps.  This class of uniformly continuous
functions is easier to deal with because under various circumstances
they possess derivatives, and their derivatives are used to reduce
the Lipschitz problem to a linear one.

We thus study, in the first two sections, the appropriate
notions of differentiability and use these to obtain some
strong results about the Lipschitz classification.

In the third section we prove Lindenstrauss' theorem on the
linearization of Lipschitz retractions to yield linear projections
(under suitable conditions).  This is not always possible and
we also show, following Lindenstrauss, that there is a Lipschitz
retraction from $\ell_\infty$ onto $c_0$.

In \S 4 we present the main result of Heinrich and Mankiewicz 
\cite{HM} which combines the techniques of \S 1 and \S 3 to
show that ``nice'' separable spaces can be Lipschitz homeomorphic
iff they are linearly isomorphic.  This is not so for general
spaces, however, and we also present the example of 
Aharoni-Lindenstrauss \cite{AL1} of non-isomorphic Lipschitz
homeomorphic Banach spaces.

In the last two sections we study general uniform homeomorphisms.
Here the situation is much more delicate and difficult,
because uniformly continous functions need not have derivatives.
It turns out that uniformly homeomorphic spaces must
have the same finite dimensional spaces.  This result is due
to Ribe \cite{Ri1}, \cite{Ri2}, and we present the simpler proof
of Heinrich and Mankiewicz \cite{HM}.  We also give in this
section a proof of an unpublished result of Enflo, that $\ell_1$
and $L_1$ are not uniformly homeomorphic.

The last section is devoted to a presentation of the recent
examples of Ribe \cite{Ri3} and Aharoni-Lindenstrauss \cite{AL2}
of uniformly homeomorphic and non-isomorphic uniformly convex,
separable spaces.
\bigskip

\head  \S1.  Existance of derivatives.\endhead 

Let $f$ be a mapping from $X$ to $Y$, and $x_0 \in X$.  We
say that $f$ is differentiable at $x_0$ if the limit
${{\partial f}\over {\partial x}} (x_0) =\lim_{t\to \infty}
{\frac 1t} (f(x_0+tx)-f(x_0))$
exists for each $x\in X$, and is linear in $x$.  (This is
usually called Gateaux differentiability, but we shall 
just call $f$ differentiable.)  The linear map $x\to {\frac
{\partial f}{\partial x}} (x_0)$ will be denoted by $(Df)_
{x_0}(x)$.

Note that if $f$ is differentiable at $x_0$ and satisfies a 
Lipschitz condition with constant $K$, then $(Df)_{x_0}$
is a bounded linear operator with $\Vert (Df)_{x_0}\Vert
\leq K$.  Moreover, if $f$ also satisfies $\Vert f(x)-f(y)
\Vert \geq {\frac 1K} \Vert x-y\Vert$ for all $x,y \in X$
then $\Vert (Df)_{x_0}(x)\Vert \geq {\frac 1K}\Vert x\Vert$,
and $(Df)_{x_0}$ is an into isomorphism.  We thus proved

\proclaim{Lemma 1.1}  If $f\ :\ X\to Y$ is a Lipschitz
embedding and $f$ is differentiable at $x_0\in X$, then
$(Df)_{x_0}$ is a linear into isomorphism.\endproclaim

This simple Lemma motivates the study of differentiability
of Lipschitz mappings, and the main result of this section
is the following theorem, proved independently by N.~Aronszajn
\cite{Ar}, J.P.R.~Christensen \cite{C}, and
P.~Mankiewicz \cite{Man}. 

\proclaim{Theorem 1.2}  A Lipschitz function from a separable
Banach space $X$ to a Banach space $Y$ with the $RNP$ is
differentiable at least at one point.\endproclaim

The requirement that $Y$ has the $RNP$ is essential for this
type of result.  Indeed, if $Y$ fails the $RNP$ there is a
Lipschitz function $f\ :\ \IR \to Y$ which is nowhere 
differentiable.

Many results on the impossibility of certain Lipschitz
embeddings follow immediately from Lemma 1.1 and Theorem 1.2,
and the known impossibility of linear embedding.  The following
examples were known before Theorem 1.2 was proved 
(Lindenstrauss \cite{L} and Enflo \cite{E1} and \cite{E3}).
But their original proofs were quite difficult and used special
properties of the spaces involved.
\smallskip
\item{(i)}  If $p<\infty$ and $q\not= p,\ \ell_q$ does not
Lipschitz embed into $\ell_p$.
\item{(ii)} Unless $p=q=2$, or $p=\infty,\ L_q$ does not
Lipschitz embed into $\ell_p$.
\item{(iii)}  If $p<\infty,\ L_q$ does not Lipschitz embed
into $L_p$ unless $2\geq q\geq p\geq 1$ or $p=q$.
\smallskip
(The case $q>2,\ p=1$ does not follow from the Theorem 
because $L_1$ does not have the $RNP$.  But it is still
true that $L_q,\ q>2$ does not Lipschitz embed into $L_1$.
See the Example at the end of \S 2.)
\smallskip
\item{(iv)}  If $X$ Lipschitz embeds into a Hilbert space, it
must be isomorphic to a Hilbert space.  Indeed, Lemma 1.1 and
Theorem 1.2 apply directly when $X$ is separable, and a space
is isomorphic to a Hilbert space if all its separable subspaces
are.
\smallskip

In order to prove Theorem 1.2 we first need to study the
differentiability of Lipschitz functions from $\IR$ or $\IR^n$
into a space $Y$ with the $RNP$.  This involves two classical 
results of Gelfand and of Rademacher.  The first proposition
was proved by I.N.~Gelfand \cite{G} when $Y$ is a separable
dual.

\proclaim{Proposition 1.3}  A Lipschitz function $f$ from $\IR$
into a space $Y$ with the $RNP$ is differentiable almost 
everywhere.\endproclaim

\demo{Proof}
Consider the set function $\mu$ defined for intervals
$[a,b]$ by $\mu([a,b])=f(b)-f(a)$.  As $f$ is a Lipschitz function,
hence absolutely continuous and of bounded variation on each
finite interval, $\mu$ can be extended in the standard way to a
$\sigma$-additive Borel measure on $\IR$ which is absolutely
continuous with respect to the Lebesgue measure.  As $Y$ 
satisfies the $RNP$, there is function $g\ :\ \IR \to Y$,
Bochner integrable on each finite interval so that $f(t)=f(0)+
\int_0^t g(s)\, ds$.

Once this representation is given, the proof that $f$ is differentiable
almost everywhere and that $f^\prime (t)=g(t)$ for almost all
$t$ is the same as in the analogous situation for real valued
functions.
\enddemo 

The second proposition is a generalization of a famous theorem
of H.~Rade-macher \cite{Ra}.  Rademacher proved it for $Y=\IR^m$,
but using Proposition 1.3 the generalization is only formal.

\proclaim{Proposition 1.4} A Lipschitz function $f$ from 
$\IR^n$ into a space $Y$ with the $RNP$ is differentiable
almost everywhere.\endproclaim 

\demo{Proof}  Recall that for a fixed $x\in \IR^n$, we denote by 
${\frac {\partial f}{\partial x}} (y)$ the derivative of
$f$ at $y$ in the $x$ direction, {\it i.e.},
$${{\partial f}\over {\partial x}}(y) =\lim_{t\to 0} (f(y+
tx)-f(y))/t$$
provided the limit exists.

By Proposition 1.3  this limit exists a.e. on each line
parallel to $x$, thus by Fubini's Theorem it exists for
almost all $y\in \IR^n$.  

Let $G$ be a countable dense additive subgroup of $\IR^n$.
By the above, the directional derivative ${\frac {\partial f}
{\partial g}} (y)$ exists 
for almost all $y\in \IR^n$ and 
for all $g\in G$, and we shall first show that except for
a set of $y$'s of measure zero it is linear in $g\in G$.

To this end let $\psi \ :\ \IR^n \to \IR$ be $C^\infty$
with compact support, and consider $f\, *\, \psi$.  On the one
hand ${\frac \partial{\partial g}} (f\, *\, \psi) = f\, *\, 
{\frac {\partial \psi}{\partial g}}$ which is linear in $g$
because $\psi $ is $C^\infty$.

On the other hand, by a simple change of variable
$${\partial \over {\partial g}} (f\, *\, \psi)(y) =\lim_{t\to
\infty} \int {{f(x+tg)-f(x)} \over t} \ \psi (y-x)\, dx$$ and
as the functions $(f(x+tg)-f(x))/t$ are bounded and converge
a.e. to ${\frac {\partial f}{\partial g}}(x)$ this limit is
also equal to ${\frac {\partial f}{\partial g}}\, *\, \psi$,
hence, the latter  is also linear in $g$, {\it i.e.}, given
any $g,h\in G$, \ $({\frac {\partial f}{\partial g}} +
{\frac {\partial f}{\partial h}} - {\frac {\partial f}{\partial
(g+h)}}) \ *\ \psi \equiv 0$.

As $\psi$ is arbitrary and $G$ is countable it follows that
for almost all $x\in \IR^n$ we have
$${{\partial f}\over {\partial g}} (x) + {{\partial f}\over 
{\partial h}} (x) = {{\partial f}\over {\partial (g+h)}} (x)
\qquad \hbox{for all}\qquad g,h\in G \ .\leqno(*)$$

Fix any $x$ for which $(*)$ holds.  We shall show that in fact,
it holds for all $g,h\in \IR^n$ and not only in $G$,\ {\it i.e.},
\ $f$ is differentiable at $x$.

For $t>0$, put $f_t(y)=(f(x+ty)-f(x))/t$.  The family $(f_t)_
{t>0}$ is equi-continuous, (in fact, they all satisfy the same
Lipschitz condition that $f$ does).  By the Arzella-Ascolli 
Theorem this set is a relatively compact set of continuous
functions.  As it has a unique limit, when $t\to 0$, on the dense subset $G$
of $\IR^n$, the limit exists for all $y\in \IR^n$, and is
continuous in $y$.  As this continuous limit is linear on
the dense subset $G$, it is linear everywhere.  This completes
the proof.
\enddemo 

The proof of Theorem 1.2 will yield more than claimed.  Not
only do Lipschitz functions have derivatives at some points,
but in fact, they are differentiable ``almost everywhere.''
As an infinite dimensional space does not carry a ``standard''
measure one needs to define what is meant by ``almost everywhere.''
Indeed, the main point in what N.~Aronszajn, J.R.P.~Christensen
and P.~Mankiewicz do is in the introduction of a useful notion
of ``almost everywhere,'' or equivalently of the complementary
sets which we shall call ``zero sets.''

Different authors introduced different notions of zero sets,
see \cite{A}, \cite{M}, \cite{C}, \cite{P}, and we shall
present the one introduced by J.P.R.~Christensen \cite{C},
but they could all be used to prove Theorem 1.2, as they all 
satisfy the following definition.

\proclaim{Definition}  Let $X$ be a Banach space.  A family
${\Cal N}$ of Borel subsets of $X$ is called a null family,
and its members are called zero sets, if:
\smallskip
\item{(i)} ${\Cal N}$ is closed under countable unions.  
\item{(ii)}
Let $A$ be a Borel set. If there exists a finite dimensional subspace $Y$ 
of $X$ so that $(A+x) \cap Y$ is a subset of $Y$ of Lebesgue measure zero for 
every $x$ in $X$, then $A$ is in ${\Cal N}$. 
\item{(iii)}  Sets in ${\Cal N}$ have empty interior.
\item{(iv)}  If $Y$ is a finite dimensional subspace of $X$
then $A\in {\Cal N}$ iff for all $x\in X$,\ $(A+x) \cap Y$ is
a subset of $Y$ of Lebesgue measure zero.
\smallskip \endproclaim

Note that when $X$ is a finite dimensional space, the family
of Borel subsets of $X$ of Lebesgue measure zero is a
null family.  In this case (iv) is a ``mild'' form of
the Fubini Theorem.  Notice that it is exactly this form of
Fubini's Theorem that we used in deducing Proposition 1.4
from the one dimensional case, Proposition 1.3.

To emphasize the fact that only (i)-(iv) are needed to prove
Theorem 1.2 and not the particular form of ${\Cal N}$ we now
complete the proof of the Theorem under the assumption that
such a null family ${\Cal N}$ exists on $X$:

\proclaim{Proposition 1.5}  Let $X$ be a separable Banach
space and let ${\Cal N}$ be a null family on $X$.  Let $Y$
satisfy the $RNP$ and assume $f\ :\ X\to Y$ is a
Lipschitz function.  Then $f$ is differentiable ``almost
everywhere'', i.e., the set of points where $f$ is not
differentiable belongs to ${\Cal N}$.\endproclaim

\demo{Proof}  Let $X_1 \subset X_2 \subset X_3 \cdots$\ be finite
dimensional subspaces of $X$ so that $X=\overline {\cup \,
X_n}$.  Define
$$D_n = \left\{ \eqalign{&x\in X\ :\ {{\partial f} \over 
{\partial a}} (x) \quad \hbox{exists for all}\cr
&a\in X_n,\quad \hbox{and is linear in}\quad a\in X_n\cr}
\right\} \ .$$

We claim that $X \backslash D_n \in {\Cal N}$.  By (iv) it
suffices to show that for each fixed $z\in X,\ ((X\backslash
D_n) + z)\cap X_n$ is a set of Lebesgue measure zero in $X_n$.
But this set is exactly the set of points in $X_n$ where
$g_z(x)=f(x+z)$ is non-differentiable as a function from 
$X_n$ to $Y$.  Thus Proposition 1.4 implies that it has zero
Lebesgue measure.

By (i) \ $\cup \, (X\backslash D_n)$ \ is also in ${\Cal N}$, and
we shall show that $f$ is differentiable at each $x\in \cap \,
D_n$,\ {\it i.e.}, in the complement of $\cup (X\backslash D_n)$.
If $x\in \cap \, D_n,\ {\frac {\partial f}{\partial a}}(x)$
exists and is linear for $a$'s in $\cup \, X_n$.  As $\cup \, X_n$
is dense in $X$, the same argument used at the end of the proof of
Proposition 1.4 shows that in fact $(Df)_x$ exists.
\enddemo 

It remains to construct a null family ${\Cal N}$, and we shall
follow \cite{C}.

\proclaim{Definition}  A Borel subset $A$ of $X$ is called a zero
set, if there is a probability measure $\mu$ on $X$ (which we call
a test measure for $A$), so that $\chi_A \, *\, \mu \equiv 0$,\ 
{\it i.e.,} \ $\int \chi_A(x-y)\, d\mu (y) =0$ for all $x\in X$.\endproclaim 

Note that this generalizes the notion of a zero set in $\IR^n$.  
Indeed, denote the Lebesgue measure on $\IR^n$ by $\lambda$.
If $\chi_A \, *\, \mu \equiv 0$ for some Borel set $A$ and a
probability measure $\mu$, integrating the equality with respect
to $\lambda$ we obtain
$$\eqalign{0&= \int\!\!\int \chi_A(x-y)\, d\mu (y)\, d\lambda (x) \cr
&= \int \biggl( \int \chi_A (x-y)\, d\lambda (x)\biggr) \, d\mu (y) \cr
&= \int \lambda (A)\, d\mu (y) = \lambda (A) \ .\cr}$$

Conversely, if $\lambda (A)=0$ take $\mu$ to be $\lambda$ 
restricted to the unit cube $C$.  Then
$$\eqalign{\chi_A \, *\, \mu &= \int \chi_A (x-y)\, d\mu (y) \cr
&= \int_C \chi_A (x-y)\, d\lambda (y) \cr
&= \lambda ((x-A) \cap C) =0 \ .\cr}$$

We now have

\proclaim{Proposition 1.6}  Let $X$ be a separable Banach space.
Then the family of zero sets defined above is a null family.\endproclaim

\demo{Proof} To see the idea of the proof that ${\Cal N}$ satisfies
(i) let us first check that if $A,B\in {\Cal N}$ so does $A\cup B$.
Indeed, let $\mu ,\nu$ be test measures for $A,B$ respectively,
and let $\eta = \mu \, *\, \nu$.  Then $\chi_A \, *\, \eta =
(\chi_A \, *\, \mu )\, *\, \nu \equiv 0$ and similarly $\chi_B \, *\,
\eta \equiv 0$.  Hence $\chi_{A\cup B} \, *\, \eta \leq (\chi_A
+\chi_B)\, *\, \eta \equiv 0$.

To carry this proof so that it applies to countable unions we need
to form infinite convolutions, and this is where the completeness
of $X$ is used.  We first need some preparations.

Let ${\Cal P}$ be the set of all Borel probability measures on $X$
with the $\omega^*$-topology induced by the bounded continuous
functions on $X$.  It is well known  that
${\Cal P}$ is separable and that its topology is given by a 
complete metric which is translation invariant, {\it i.e.},\ 
$d(\mu \, *\, \eta, \nu \, *\, \eta)= d(\mu ,\nu)$ for all
$\mu ,\nu , \eta \in {\Cal P}$. (See {\it e.g.}, \cite{Bi}.)

Note that if $A \in {\Cal N}$ and $\mu$ is a test measure for
$A$, then so is every translate of $\mu$, and so is the restriction
${\frac 1{\mu (B)}} \mu \vert_B$ for any Borel subset $B$ so that
$\mu (B) >0$.  Thus by appropriate translation and restriction we
can find a test measure for $A$ which is supported in an 
arbitrarily small ball centered at 0, {\it i.e.}, a test
measure for which $d(\mu, \delta_0)$ is arbitrarily small.

Assume now that $A_n\in {\Cal N}$ and let $\mu_n$ be test
measures for $A_n$ with $d(\mu_n , \delta_0) < 2^{-n}$.  Then
$\mu_1 \, * \ldots *\, \mu_n$ is a $d$-Cauchy sequence in
${\Cal P}$ and let $\mu$ be its limit.  For any $n$, \ 
$\mu =\mu_n \, *\, \eta_n$ where $\eta_n$ is the infinite
convolution of all the $\mu_k$'s except the $n^{th}$.
Thus $\chi_{A_n} \, *\ \mu = (\chi_{A_n} \, *\, \mu_n)\, * \,
\eta_n \equiv 0$, and hence also $\chi_{\cup A_n}\ *\ \mu
\leq \sum \chi_{A_n} \ *\ \mu \equiv 0$, and $\cup A_n \in
{\Cal N}$.  This proves (i), and (ii) is obvious.  To prove
(iii) assume $A\in {\Cal N}$ has a non-empty interior.  We
can cover $X$ by a countable number of translates of $A$,
and by (i) and (ii) this will yield that $X\in {\Cal N}$.  
But this is impossible as $\chi_X \ *\ \mu \equiv \chi_X$
for all probability measures $\mu$.

It remains to check (iv).  Let $\mu$ be a probability measure
supported in $Y$ and equivalent to the Lebesgue measure on 
$Y$, and we shall check that it is a test measure for $A$.
Indeed
$$\eqalign{\chi_A \ *\ \mu &= \int \chi_A (x-y)\, d\mu (y) = 
\mu (x-A) \cr
&= \mu ((x-A) \cap Y) = 0 \ .\cr}$$
This completes the proof of the Proposition, hence also of
Theorem 1.2.
\enddemo 

\demo{Remark}  The methods of this section seem to
yield only into linear isomorphisms even if the original
Lipschitz map was an onto Lipschitz homeomorphism.  If we 
knew that there is a point $x_0$ where $f$ is differentiable
and so that, at the same time, 
$f^{-1}$ is differentiable at $f(x_0),\ (Df)_{x_0}$
would have been a linear onto isomorphism with inverse $(Df^{-1})_
{f^{-1}(x_0)}$.  This led Christensen to raise the following
problem:  Assume $X$ and $Y$ are separable spaces and that 
$f\ :\ X\to Y$ is a Lipschitz homeomorphism.  Does $f$
carry null sets in $X$ to null sets in $Y$?  If the answer
is positive, the results of this section imply that ``most''
points $x_0$ in $X$ will be such that both $f$ and $f^{-1}$ are
differentiable in $x_0$ and $f(x_0)$ respectively.
\enddemo 

\head  \S 2.  $\omega^*$-derivatives.\endhead 

It follows from Theorem 1.2 that Lipschitz mappings into 
separable dual spaces have derivatives, and since the
derivative of a Lipschitz embedding is a linear isomorphism,
it follows that when a separable space Lipschitz embeds into
a separable dual $Y$, it embeds linearly into it.

In this section we show that this last result is true even
when $Y$ is a non-separable dual space.  As $Y$ will
no longer have to have the $RNP$, the Lipschitz embedding will not
have to be differentiable anywhere.  We shall define a weaker
notion of a derivative --- the $\omega^*$-derivative, and it
will turn out that the embedding will have $\omega^*$-derivative 
``almost everywhere''.  The $\omega^*$-derivative will be
a bounded linear operator, but it will not be automatically
bounded from below.  It turns out, however, that it will be
bounded from below ``almost everywhere''.  The results of
this section are taken from S.~Heinrich and P.~Mankiewicz
\cite{HM}.

\proclaim{Definition}  Let $Y=Z^*$ be a dual space and
$f\ :\ X\to Y$.  We say that $f$ is $\omega^*$-differentiable
at $x_0\in X$ if the $\omega^*$-limit 
$$(D^*f)_{x_0}(x) = \omega^* - \lim_{t\to 0} (f(x_0+tx) -f(x))
/t$$
exists for all $x\in X$, and is linear in $x$.  The linear 
mapping $(D^*f)_{x_0}$ is called the $\omega^*$-derivative
of $f$ at $x_0$.\endproclaim

If $f$ satisfies a Lipschitz condition with constant $K$,
the $\omega^*$-lower-semi-continuity of the norm implies
that $(D^*f)_{x_0}$ is a linear operator of norm at most $K$.

\proclaim{Theorem 2.1}  A Lipschitz map $f$ from a separable
Banach space into the dual $Y=Z^*$ of a separable space $Z$,
has a $\omega^*$-derivative almost everywhere.  Moreover,
if $f$ is a Lipschitz embedding, $(D^*f)_{x_0}$ is an
isomorphism for almost all $x_0$.\endproclaim

\demo{Proof}  By ``almost everywhere'' we mean, of course,  
``except for a zero set'' in the sense described in the end of 
\S 1.  As $Z$ is separable, most of the proof is routine
extension of the ideas used in \S 1, applied to the scalar
Lipschitz mappings $\varphi_n (x) = \langle z_n, f(x) \rangle$
where $(z_n)$ is a dense subset of the unit ball of $Z$.
The only new ingredient is in the proof of the last part of
the Theorem.  Here also the whole idea can be seen in the 
one-dimensional case $X=\IR$.  The extension to general
separable $X$ follows the same route taken in \S 1.  Hence,
we treat this case only.  So assume $f\ :\ \IR \to Z^*$ 
satisfies a Lipschitz condition with constant $K$, and
without loss of generality normalize so that $\Vert f(t) -
f(s)\Vert \geq \vert t-s\vert$ for all $t,s\in \IR$.  To
simplify the notation write $f^\prime (t_0)$ for $(D^*f)_
{t_0}$ when the latter exists.  Thus, in this case
$$\omega^* - \lim_{\lambda \to 0} (f(t_0 +\lambda t) -
f(t_0))/\lambda = t f^\prime (t_0)\ .$$
By the earlier parts of the Theorem we already know that
$f^\prime (t)$ exists almost everywhere, and then, of
course, $\Vert f^\prime (t) \Vert \leq K$.  We need to
show that $\Vert f^\prime (t)\Vert \geq 1$ for almost
all $t$.  

If this were false, find $\delta <1$ so that $A=\{ t\ :\ 
\Vert f^\prime (t)\Vert \leq \delta \}$ has positive
measure, and let $t_0$ be a density point of $A$.  Fix
$\epsilon >0$ so small that $$m([t_0, t_0+\epsilon] \cap A)
\geq \epsilon ( 1- {\frac {1-\delta}{3K}} )$$ and let
$B=[ t_0,t_0 +\epsilon ] \backslash A$.  Then $m(B) \leq
\epsilon \ {\frac {1-\delta}{3K}}$ and $\Vert f^\prime (t)
\Vert \leq K$ on $B$. 

As $\Vert f(t_0 +\epsilon )-f(t_0)\Vert \geq \epsilon$, find
$z\in Z,\ \Vert z\Vert =1$ with $\langle z,f(t_0+\epsilon)
-f(t_0)\rangle \geq \epsilon (1+\delta )/2$, and define
$\varphi (t) = \langle z,f(t)\rangle$.  $\varphi$ is a
scalar Lipschitz function, hence the integral of its derivative.
Thus
$$\eqalign{\epsilon (1+\delta )/2 &\leq \varphi (t_0 +\epsilon )
-\varphi (t_0) = \int\limits_{t_0}^{t_0 +\epsilon} \varphi^\prime (t)\, dt \cr
&\leq \int_B \Vert \varphi^\prime (t)\Vert +\int_A
\Vert \varphi^\prime (t)\Vert \cr
&\leq K\, m(B) + \epsilon \delta \leq \epsilon \bigl( {{1-\delta}\over 3}
+ \delta \bigr) \cr}$$
a contradiction.
\enddemo 

Theorem 2.1 applies to dual of separable spaces.  But for our
purposes the separability requirement can be overcome by the
following lemma which is of interest for its own sake.

\proclaim{Lemma 2.2}  Let $E$ be a separable subspace of the dual
space $F^*$.  Then there is a subspace $Y$ of $F^*$, containing
$E$ which is isomorphic to the dual of a separable space $Z$.\endproclaim

For the proof, see \cite{HM}, Proposition 3.4.
\medskip

Combining this result with Theorem 2.1 we obtain

\proclaim{Corollary 2.3}  Let $X$ be a separable Banach space.
If $X$ Lipschitz embeds into a dual space, it embeds linearly
into it.  In particular, if $X$ Lipschitz embeds into any space
$Y$, it embeds linearly into $Y^{**}$.\endproclaim

\demo{Proof}  Assume $f\ :\ X\to F^*$ is a Lipschitz embedding.  As
$X$ is separable $f(X)$ is contained in a separable subspace
$E$ of $X^*$, and find $Y,Z$ as in Lemma 2.2.  Now use
Theorem 2.1 to find $x_0$, so that $(D^*f)_{x_0}$ is a
linear isomorphism into $Y$.  (Note that $(D^*f)_{x_0}$ is
taken with respect to the $\omega^*$ topology on $Y=Z^*$
and not with respect to the original $\omega^*$ topology on
$F^*$.)
\enddemo 

\demo{Example}  We can now show that $L_q$ does
not Lipschitz embed in $L_1$ for $q>2$.  Indeed, by the
Corollary this would have implied that $L_q$ embeds linearly
into the $L$-space $L_1^{**}$ which is impossible.
\enddemo 

\head  \S 3.  Linearization of Lipschitz retractions.\endhead 

The main result of this section is the following Theorem of
J. Lindenstrauss \cite{L}.

\proclaim{Theorem 3.1}  Let $X$ be a subspace of $Y$ so that
there is a Lipschitz retraction from $Y$ onto $X$.  If $X$ is
complemented in $X^{**}$ it is complemented in $Y$.\endproclaim

Before we give the detailed proof let us explain the idea.
Assume $f\ :\ Y \buildrel \hbox{onto} \over \longrightarrow X$
is the Lipschitz retraction, and that $f$ is differentiable
at a point $x_0 \ \underline{\hbox{in } X}$.  Then $(Df)_{x_0}$ is
a linear operator from $Y$ into $X$ which is a projection.
Indeed, $f$ is the identity on $X$, thus for every $x\in X$\ 
and every $t$,
$(f(x_0+tx)-f(x_0))/t =x$, hence also $(Df)_{x_0}(x)=x$.

Our assumptions do not guarantee that $f$ is differentiable
anywhere --- and certainly not in a point of the ``small''
subset $X$ of $Y$.  We thus need a ``smoothing'' procedure.
This will be done on each finite dimensional subspace of 
$Y$ (see Lemma 3.3), and then we shall use a routine compactness
argument to obtain an operator $T\ :\ Y\to X^{**}$ so that
$T$ is the identity on $X$.  As $X$ is complemented in $X^{**}$,
say by a projection $P,\ PT\ :\ Y\to X$ will be the desired
projection.

>From this description it might appear that the assumption that
$X$ is complemented in $X^{**}$ is inessential, and is just a 
technical drawback of our method of proof.  But this is not the
case as the following example of Lindenstrauss \cite{L} shows.

\proclaim{Example}  There is a Lipschitz retraction from $\ell_
\infty$ onto $c_0$.\endproclaim
 
Given $x=(a_1,\ldots ,a_n, \ldots )\in \ell_\infty$, denote by
$d(x)$ its distance to $c_0$,\ {\it i.e.},\ $d(x) =\lim \sup
\vert a_n\vert$.  Define now $f\ :\ \ell_\infty \to c_0$ by
$$(f(x))(n) = \cases 0&\text{$\vert a_n\vert \leq d(x)$}\cr
{{a_n}\over {\vert a_n\vert}}(\vert a_n\vert -d(x))&\text{$\vert
a_n\vert > d(x)$}\endcases\ .$$
One easily checks that $f$ is a retraction from $\ell_\infty$
onto $c_0$ and that it satisfies a Lipschitz condition with
constant 2.
\medskip

The example above  is in fact a special case of a more general
theorem of Lindenstrauss \cite{L} which we just quote without
proof.

\proclaim{Theorem 3.2}  If $K$ is a compact metric space then
$C(K)$ is an absolute Lipschitz retract, {\it i.e.}, if $Y$ is
any metric space containing $C(K)$ there is a Lipschitz retraction
from $Y$ onto $C(K)$.\endproclaim

We now pass to the proof of Theorem 3.1.  The main step in the
proof is contained in the following Lemma, where Lipschitz
maps on finite dimensional spaces are ``linearized''. 

\proclaim{Lemma 3.3} Let $Z$ be a finite dimensional space and
let $E$ be a subspace of $Z$.  Let $f$ be a Lipschitz function
from $Z$ into a Banach space $X$ so that $f\vert_E$ is a given 
linear operator $S$.  Then there is a linear operator $T\ :\ Z\to
X^{**}$ so that $T\vert_E=S$ and whose norm is at most the
Lipschitz constant of $f$.\endproclaim

\demo{Proof}  Write (algebraically) $Z= E\oplus Y$, and assume 
$\dim Y=m$ say.  We first claim that without loss of generality
we can assume that  the directional derivatives ${\frac {\partial f}
{\partial e}}$ in the direction of $e\in E$ exist everywhere,
are continuous in $z$, and are linear in $e\in E$.

Indeed, let $\varphi \geq 0$ be a $C^\infty$ function on $E$
with compact support so that $\int \varphi =1$ and $\varphi (x)=
\varphi (-x)$ for all $x\in E$, and replace $f$ by the function
$F$ given by
$$F(z) = \int\limits_E f(z+x) \varphi (x)\, dx\ .$$

It is easy to check (using change of variable as below) that
$F$ is differentiable as required.  To check the continuity
of the directional derivatives, fix $e\in E$ and $t$.  By a
change of variable write
$$F(z+te) = \int\limits_E f(z+x) \varphi (x-te)\, dx \ .$$

Thus for $z,z^\prime \in Z$ we have 
$$\eqalign{\Vert (F(z+te)-F(z))/ t-(F &(z^\prime +te)-
F(z^\prime))/ t \Vert =\cr &= \Vert \int\limits_E (f(z+x)-
f(z^\prime +x))(\varphi (x-te)-\varphi (x)) / t\, dx \Vert \cr
&\leq K\Vert z-z^\prime \Vert \ \int_E \ \Vert (\varphi (x-te) -
\varphi (x))/t \Vert \, dx \cr}$$
(where $K$ is the Lipschitz constant of $f$).  Passing to the limit
as $t\to 0$ gives
$$\Vert {{\partial F}\over {\partial e}} (z) - {{\partial F}\over
{\partial e}} (z^\prime) \Vert \leq K\Vert z-z^\prime \Vert \ \int_E
\ \Vert {{\partial \varphi}\over {\partial e}} (x) \Vert \, dx$$
hence $\frac {\partial F}{\partial e}$ satisfies a Lipschitz
condition.

To see that the restriction of $F$ to $E$ is $S$, we use the
linearity of $f\vert_E$:  If $e\in E$,
$$\eqalign{F(e) &= \int\limits_E f(e+x)\varphi (x)\, dx \cr
&= \int_E (f(e)+f(x)) \varphi (x)\, dx = f(e) \cr}$$
because $\int \varphi (x)=1$ and $\int f(x)\varphi (x)\, dx =0$
because $f$ is linear on $E$, so $f(x)=-f(-x)$ while $\varphi
(x)=\varphi (-x)$.

>From now on we can thus assume $f$ is already ``smooth'' in the
$E$ directions.

To prove the Lemma we shall now use a ``smoothing kernel'' in
the $Y$ direction.  Let $\psi \geq 0$ be a $C^\infty$ function on
$Y$ with compact support so that $\int_Y \psi =1$.  Define
$$f_n(z) = n^m \int\limits_Y f(z+y) \psi (ny)\, dy\ .$$
As $f$ is already smooth in the $E$ directions, each $f_n$ is
differentiable so let $T_n=(Df_n)_0$ be its differential at 0.
The sequence of operators $(T_n)$ is uniformly bounded from $Z$
to $X\subset X^{**}$ so fixing a free ultra-filter $U$ on $\text{\bf N}$
define $Tz = \omega^* -\lim_U T_nz$.  We need only check that
$T\vert_E =S$. So fix $e\in E$, and then $Te =\lim T_ne = \lim
{\frac {\partial f_n}{\partial e}} (0)$.  But
$${\frac {\partial f_n}{\partial e}} (0) = n^m \int_Y
{\frac {\partial f}{\partial e}} (y) \psi (ny)\, dy =
\int_Y {\frac {\partial f}{\partial e}} (n^{-1}y) \psi
(y)\, dy \ ,$$
and, of course, ${\frac {\partial f}{\partial e}} (0) = Se$
because $f\vert_E=S$.  Thus $$Te-Se=\lim \int_Y (
{\frac {\partial f}{\partial e}} (n^{-1} y) -
{\frac {\partial f}{\partial e}} (o)) \psi (y)\, dy$$ which is
zero because $\psi$ has a compact support and $\frac {\partial f}
{\partial e}$ is continuous.
\enddemo 

\demo{Proof of Theorem 3.1}  The proof now follows
easily from the Lemma.  Let $Y_\alpha$ be  a net, directed by
inclusion, of finite dimensional subspaces of $Y$ so that $Y=\cup Y_\alpha$.
For each $\alpha$ use the Lemma with $Z=Y_\alpha$ and $E=E_\alpha =
Y_\alpha \cap X \subset Z$.  Note that since $f$ is a retraction
of $Y$ onto $X,\ f\vert_{E_\alpha}$ is the identity on $E_\alpha$.
By the Lemma, find a linear operator $T_\alpha \ :\ Y_\alpha \to
X^{**}$ so that $T_\alpha \vert_{X\cap Y_\alpha}$ is the identity
on this space, and so that $\Vert T_\alpha \Vert \leq K$, the
Lipschitz constant of $f$.  By compactness $(T_\alpha)$ has a
$\omega^*$-convergent subnet, thus its limit $T$ is an operator
from $Y$ to $X^{**}$ with $T\vert_X =id_X$.  Now take $PT\ :\ 
Y\to X$ as the desired projection, where $P\ :\ X^{**}\to X$ is
the given bounded linear projection.
\enddemo 

\head  \S 4.  Linear isomorphism between Lipschitz
equivalent spaces.\endhead 

The differentiation theory of \S 1 yields that if two ``nice''
spaces are Lipschitz equivalent, they embed linearly into each
other, but we couldn't prove that they are actually isomorphic
(see the Remark at the end of \S 1).  It turns out, however,
that one can combine the results of \S 1 with the retraction
linearization of \S 3 to obtain that the ``nice'' spaces 
embed linearly as {\it complemented} subspaces of each other.
Thus, we obtain linear isomorphism results for a large class
of spaces for which the ``decomposition scheme'' holds.  After
presenting this method, which is due to S.~Heinrich and P.~
Mankiewicz \cite{HM}, we present the example of I.~Aharoni
and J.~Lindenstrauss \cite{AL1} of two (``non-nice'') spaces
which are Lipschitz equivalent but not linearly isomorphic.

\proclaim{Theorem 4.1}  Let $f$ be a Lipschitz embedding of
$X$ into a space $Y$ so that there is a Lipschitz retraction
from $Y$ onto the image $f(X)$.  Assume that $f$ is differentiable
at a point $x_0\in X$ and that $X$ is a Lipschitz retract of 
$X^{**}$.  Then there is a Lipschitz retraction from $Y$ onto
$(Df)_{x_0}(X)$.  In particular, if $X$ is linearly complemented
in $X^{**}$, there is a linear projection from $Y$ onto
$(Df)_{x_0}(X)$.\endproclaim

\demo{Proof}  The last claim follows from the first and Theorem 3.1,
 so we need only prove the first. WLOG $x_0=0$ and $f(0)=0$,
and denote $(Df)_0$ by $D$.  Let $\pi \ :\ Y\to f(X)$ be the
Lipschitz retraction and define $g$ and $g_n$ from $Y$
into $X$ by $g=f^{-1} \circ \pi$ and $g_n(y) =ng(y/n)$.  The
functions $g_n\ :\ Y\to X\subset X^{**}$ are uniformly 
Lipschitz and for each $y\in Y$ the sequence $g_n(y)$ is bounded.
We can thus find a $\omega^*$-limit point $h\ :\ Y\to X^{**}$ of
the $g_n$'s which, by the $\omega^*$-lower semi-continuity
of the norm is again Lipschitz with the same constant.

We claim that $h(Dx)=x$ for all $x\in X$.  Indeed, put $y=Dx$, 
{\it i.e.},\ $y=\lim_n nf(x/n)$.  As the $g_n$'s are uniformly
Lipschitz, we deduce that $$g_n(y)-g_n(nf(x/n)) \to 0\ .$$  But 
$g_n(y)$ has a subnet converging $\omega^*$ to $h(y)$ while
$$g_n(nf(x/n))=ng(f(x/n))=x\ ,$$
{\it i.e.},\ $x=h(y)=h(Dx)$.
The desired retraction is now $D \circ \rho \circ h$, where $\rho$
is the retraction from $X^{**}$ onto $X$. 
\enddemo 

Next we formulate one isomorphism result that follows from
the Theorem.  We refer the reader to \cite{HM} for many more
variations and refinements of the same theme.

\proclaim{Corollary 4.2}  Let $X$ and $Y$ be separable reflexive
spaces isomorphic to their squares.  If $X$ and $Y$ are Lipschitz
equivalent, they are linearly isomorphic. \endproclaim

\demo{Proof}  By Theorem 1.2 all the conditions of Theorem 4.1 are satisfied.
Thus $X$ and $Y$ are each linearly isomorphic to a complemented
subspace of the other, and the result follows from Pe\l czy\`nski's
decomposition method.
\enddemo 

As the next example shows, some conditions are necessary to deduce
linear isomorphism from Lipschitz equivalence.  But the spaces in
the example are not ``nice''.  They are non-separable and do not
have the Radon-Nykodim property.  In particular it is unknown if
separable examples exist or if reflexive examples exist.

\proclaim{Example}  {\rm (\cite{AL1}).}  There are Lipschitz
equivalent spaces $X$ and $Y$ which are not linearly isomorphic.
In fact $Y$ does not embed linearly into $X$. \endproclaim

Let $\{ N_\gamma \ :\ \gamma \in \Gamma \}$ be an uncountable
collection of infinite subsets of natural numbers so that
$N_\gamma \cap N_\beta$ is finite for all $\gamma \not= \beta$.
Take $X$ to be the subspace of $\ell_\infty$ spanned by $c_0$
and the characteristic functions $\chi_{N_\gamma}$ of the sets
$N_\gamma$.  Let $Y$ be $c_0(\Gamma)$.  The coordinate
functionals on $\ell_\infty$ are a countable family in 
$X^*$ which separates the points of $X$.  As no countable
family in $Y^*$ separates the points of $Y$, it cannot be 
linearly embedded into $X$.

The quotient space $X/c_0$ is isometric to $c_0(\Gamma)$.
Indeed, for any $\gamma_1,\ldots ,\gamma_n \in \Gamma$ and
$a_1,\ldots ,a_n$, 
$$\Vert \sum_1^n a_j \chi_{N_{\gamma_j}}
\Vert_{X/c_0} = \max \vert a_j\vert \ ,$$ 
{\it i.e.}, the images of $\chi_{N_\gamma}$ in the quotient are
isometrically equivalent to the unit vectors of $c_0(\Gamma)$.
Let $q\ :\ X\to X/c_0$ be the quotient map.  The heart of the
construction is to show that $q$ admits a Lipschitz ``lifting'',
\ {\it i.e.}, there is a Lipschitz map $f\ :\ X/c_0 \to X$
so that $q \circ f$ is the identity on $X/c_0$.  Once this
is shown it follows easily that $X$ and $Y$ are Lipschitz
equivalent.  Indeed, $X$ is Lipschitz equivalent to $c_0 \oplus
X/c_0$ via the map $x\to (x-f(q(x)), q(x))$ whose inverse is
$(y,z)\to y+f(z)$.  But $c_0 \oplus X/c_0 = c_0\oplus c_0(\Gamma)$
is isometric to $c_0(\Gamma) =Y$.

We shall define $f$ on $c_0(\Gamma)^+$, the nonnegative elements
of $c_0(\Gamma)$ so as to satisfy $\Vert f(y)-f(t)\Vert \leq
\Vert y-z\Vert$ whenever $y,z\in c_0(\Gamma)^+$.  We shall then
define, for any $y\in c_0(\Gamma),\ f(y)=f(y^+)- f(y^-)$ where
$y=y^+ -y^-$ is the cannonical representation of $y$ as a difference
of two disjointly supported nonnegative terms.  It follows that
$$\eqalign{\Vert f(y)-f(z)\Vert &\leq \Vert f(y^+)-f(z^+)\Vert +
\Vert f(z^-)-f(y^-)\Vert \cr
&\leq 2 \max \biggl( \Vert y^+ -z^+\Vert, \Vert y^- -z^-\Vert
\biggr) \leq 2
\Vert y-z\Vert \cr}$$
and $f$ is a Lipschitz map with constant 2.
Thus fix $y\geq 0$ in $c_0(\Gamma)$ and assume $y=\sum_{j=1}^\infty
a_j e_{\gamma_j}$ where $a_1 \geq a_2 \geq  \cdots \ $. 

Let $M_1 =N_{\gamma_1}$ and define inductively $M_n =N_{\gamma_n}
\backslash \cup_{j< n} N_{\gamma_j}$.  Define now
$$f(y) = \sum a_j \chi_{M_j} \ .$$
Note that each $\chi_{M_j} \in X$ as it differs from $\chi_{N_{\gamma_j}}$
by a $c_0$ element.  Also $q(\chi_{M_j} )=q(\chi_{N_{\gamma_j}}) =
e_{\gamma_j} \in c_0 (\Gamma)$, so $q(f(y))=y$.  To see that $f$ 
satisfies a Lipschitz condition we give another formula for
$f(y)$:  Given $n,\ (f(y))_n$, the $n^{th}$ coordinate of $f(y)$,
is equal to $a_i$ iff $n\in N_{\gamma_i}$ but $n\notin N_{\gamma_1}
\cup \cdots \cup N_{\gamma_{i-1}}$.  So let $A_n$ be the closed
subspace of $c_0(\Gamma)$ given by $A_n = \overline {sp} \{ e_\gamma
\ :\ n\in N_\gamma \}$.  It follows from the above and the 
monotonicity of the $a_i$'s that if $y\in c_0 (\Gamma)^+$ we have:
$$(f(y))_n = \hbox{ dist }\{ y,A_n\} \ .$$
It is clear from this formula that given $y,z \in c_0(\Gamma)^+$, 
$$\vert  (f(y))_n -(f(z))_n \vert = \vert \hbox{ dist}\{ y,A_n\} -
\hbox{dist}\{ z,A_n \} \vert \leq \Vert y-z\Vert \ ,$$  {\it i.e.,}
\ $\Vert f(y)-f(z)\Vert \leq \Vert y-z\Vert $. 

\head \S 5.  Uniform homeomorphisms.\endhead

As uniformly continuous functions do not have derivatives in general,
their study requires bona fide metric, geometric and topological
arguments, and usually cannot be reduced to the linear theory.

The following simple lemma is usually the tool by which we get some
initial control on a uniformly continuous function.  We say that
$f\ :\ X\to Y$ is Lipschitz for large distances if for each
$\delta > 0$ there is a $K=K(\delta)$  so that $\Vert f(x)-
f(y)\Vert \leq K\Vert x-y\Vert$ for all $x,y\in X$ satisfying
$\Vert x-y\Vert \geq \delta$.

\proclaim{Lemma 5.1}  A uniformly continuous function between
Banach spaces is Lipschitz for large distances.\endproclaim

\demo{Proof}  Given $\delta$, choose first $M$ so that $\Vert f(a) -
f(b) \Vert \leq M$ for all $a,b\in X$ satisfying $\Vert a-b\Vert
< \delta$, and let $K(\delta) =2M/\delta$.  Given $x,y$ with
$\Vert x-y\Vert \geq \delta$, let $x=a_0,a_1, \ldots ,a_m =y$
be points in $X$ so that $\Vert a_{j+1} -a_j\Vert < \delta$,
and $m=[ 2\Vert x-y\Vert /\delta ]$.  Then
$$\Vert f(x) -f(y)\Vert \leq \sum_{i=1}^m \Vert f(a_i) -
f(a_{i-1})\Vert \leq mM \leq K(\delta)\Vert x-y\Vert \ .$$
\enddemo 

The next theorem (from \cite{HM}) will enable us to use some of
the tools that we developed for Lipschitz maps in the uniformly
continuous case.

\proclaim{Theorem 5.2}  If $X$ and $Y$ are uniformly homeomorphic,
they have Lipschitz equivalent ultra powers.\endproclaim 

\demo{Proof}  Let $f\ :\ X\to Y$  be a uniform homeomorphism.  By
Lemma 5.1 there is a constant $K$  so that $\Vert f(x_1)-f(x_2)
\Vert \leq K\Vert x_1 -x_2\Vert$ and $$\Vert f^{-1}(y_1) -
f^{-1}(y_2)\Vert \leq K\Vert y_1-y_2\Vert$$ for all 
$x_1,x_2\in X$ and $y_1,y_2\in Y$ satisfying $\Vert x_1-x_2
\Vert ,\Vert y_1-y_2\Vert \geq 1$.  Define $f_n(x)= {\frac 1n}
f(nx)$.  Then the modulus of continuity of $f_n$ is not worse than 
that of $f$, and $f_n$ already  satisfies the $K$-Lipschitz 
condition $\Vert f_n(x_1)-f_n(x_2)\Vert \leq K\Vert x_1 -
x_2\Vert$ whenever $\Vert x_1 -x_2\Vert \geq {\frac 1n}$.  As
$f_n^{-1}(y)= {\frac 1n} f^{-1}(ny)$ the same remarks apply
to $f_n^{-1}$.

Let $U$  be a free ultra filter on the natural numbers, and
define $F=(f_n)$ to be the natural map from $(X)_U$, the
ultra power of $X$ onto $(Y)_U$.

We shall show that $F$ satisfies a Lipschitz condition with
constant $K$.  Fix $x=(x_1,\ldots)$ and $z=(z_1,\ldots )$
in $(X)_u$, and choose, by uniform continuity of $f$, a 
$\delta >0$ so that $\Vert f(a)-f(b)\Vert \leq K\Vert x-z\Vert$
for all $a,b\in X$ satisfying $\Vert a-b\Vert < \delta$.
Then we also have $\Vert f_n(a)-f_n(b)\Vert \leq K\Vert x-z
\Vert$ for all $n$.  As there are only finitely many $n$'s
for which $\delta \leq \Vert x_n -z_n\Vert \leq {\frac 1n}$, we 
see, by the choice of $\delta$ and $K$, that 
for each $\epsilon >0$, \  $\Vert f_n(x_n)-
f_n(z_n)\Vert \leq K\Vert x-z\Vert +\epsilon$ for all but finitely 
many $n$'s, \ {\it i.e.},\ $\Vert F(x)-F(z)\Vert \leq K
\Vert x-z\Vert$.
Similar arguments hold for $F^{-1}$.
\enddemo 

We shall now deduce two results of Ribe \cite{Ri1}, \cite{Ri2}  from
Theorem 5.2.  These are only two examples of the consequences
of Theorem 5.2 and its variations, and we refer the reader to
the fundamental paper of Heinrich and Mankiewicz \cite{HM}
for many more details.

\proclaim{Theorem 5.3}  If $X$ is uniformly homeomorphic
to a ${\Cal L}_p$ space, $1< p<\infty$, then it is a ${\Cal L}_p$
space itself.\endproclaim

\demo{Proof}  $X$ is a ${\Cal L}_p$ space iff $(X)_u$ is such a space
for some ultra filter $U$, and by Theorem 5.2 there is a $U$ for 
which $(X)_u$ is Lipschitz equivalent to an ultra power of a
${\Cal L}_p$ space, {\it i.e.}, to a ${\Cal L}_p$ space.  It thus
remains to show that a space $Z$, Lipschitz equivalent to a
${\Cal L}_p$ space $Y$ is itself a ${\Cal L}_p$ space.  Let
$f\ :\  Z\to Y$ be the uniform homeomorphism.

If $Z$ is separable, this follows immediately from Theorems
1.2 and 4.1.  Indeed, $Y$ is reflexive (because $1< p<\infty$),
so $f$ has a point of differentiability, hence $Z$ is isomorphic
to a complemented subspace of $Y$.

The general case is reduced easily to the separable one, by 
showing that every separable subspace $Z_0\subset Z$ is contained in
a separable subspace $W\subset Z$ which is Lipschitz equivalent
to a ${\Cal L}_p$ space.  Indeed, define inductively sequences
$Z_0\subset Z_1 \subset \cdots$ in $Z$ and $Y_0 \subset Y_1 \subset
\cdots $ in $Y$ as follows:  Having defined $Z_n$, let $Y_n$
be a separable ${\Cal L}_p$ subspace of $Y$ containing $f(Z_n)$,
and having defined $Y_n$ let $Z_{n+1}$ be any subspace of $Z$
containing $f^{-1}(Y_n)$.  Then $\overline{\cup Y_n}$ is a 
${\Cal L}_p$ space, Lipschitz equivalent to $W={\overline {\cup Z_n}}$.
\enddemo 

\proclaim{Theorem 5.4}  If $X$ and $Y$ are uniformly homeomorphic,
there is a constant $C\geq 1$ so that for each finite dimensional 
subspace $E$ of $X$ there is a subspace $F$ of $Y$ with $d(E,F)
\leq C$.\endproclaim 

\demo{Proof}  By Theorem 5.2 there is an ultra filter $U$ and a Lipschitz
equivalence ${f : (X)_u \to (Y)_u}$.  
Fixing $E\subset X$ and considering it as a subspace of $(X)_u$,
we have by Corollary 2.3  that $E$ embeds linearly into
$((Y)_u)^{**}$.  The result now follows from local reflexivity and
the local structure of ultra powers.
\enddemo 

Theorem 5.4 gives many examples of non-uniformly homeomorphic
Banach spaces and leads naturally to the question of whether
some converse to it holds, {\it i.e.}, whether identical local
structures for two separable spaces imply that they are 
uniformly homeomorphic.  As the following example shows there
must be some additional requirements for such a converse to hold.
The example is due to P.~Enflo (unpublished) and was shown to
me by J.~Lindenstrauss.

\proclaim{Example}  $L_1$ and $\ell_1$ are not uniformly
homeomorphic.\endproclaim

We start with some observations for general Banach spaces, and then
specialize to $L_1$ and $\ell_1$.

Assume $f\ :\ X\to Y$ is a uniform homeomorphism.  As $f^{-1}$ is
uniformly continuous it is Lipschitz for large distances, and
there is a constant $L>0$ so that
$$\Vert x-y\Vert \leq \max \{ 1, L \Vert f(x)-f(y)\Vert \} \leqno(*)$$
for all $x,y\in X$.

Also, for each $\delta >0$, let $K(\delta )$ be the smallest 
Lipschitz constant of $f$ for distances above $\delta$, \ {\it i.e.},
the smallest constant $K(\delta )$ so that
$$\Vert x-y\Vert \leq \delta \Longrightarrow \Vert f(x)-f(y)
\Vert \leq K (\delta ) \Vert x-y\Vert \ . $$
(Such a $K(\delta )$ exists by Lemma 5.1.)

Obviously, $K(\delta )$ decreases as $\delta$ increases, and
let $K=\lim_{\delta \to \infty} K(\delta )$.  Then $K>0$.  In
fact, whenever $\Vert x-y\Vert \geq \delta > 1$, $(*)$ and the
definition of $K(\delta )$ imply that $\Vert f(x)-f(y)\Vert
\leq K(\delta )\Vert x-y\Vert \leq K(\delta )L\Vert f(x)-f(y)
\Vert$ so that $K\geq 1/L$.

Fix $0< \epsilon < {\frac 12}$, to be specified at the end of
the proof, and choose $\delta >1$ large enough so that
$K(\delta )\leq (1+\epsilon) K$.  Using the minimality of
$K(2\delta )$, fix $x,y\in X$ with $\Vert x-y\Vert \geq 2\delta$
so that
$$\Vert f(x)-f(y)\Vert \geq (1-\epsilon) K(2\delta )\Vert 
x-y\Vert \geq
(1-\epsilon) K\Vert x-y\Vert \ .$$

Let $U =\{ u\in X\ :\ \Vert u-x\Vert =\Vert u-y\Vert =\Vert x-y
\Vert / 2\}$ be the set of metric midpoints between $x$ and $y$,
and put
$$\eqalign{V_\epsilon =\{ v\in Y\ :\ (1-4\epsilon) \Vert f(x)-f(y)\Vert/2 
&\leq \Vert f(x)-v\Vert, \Vert f(y)-v\Vert \cr &\leq (1+4\epsilon )
\Vert f(x)-f(y)\Vert /2 \} \ .\cr}$$
The set $V_\epsilon$ is the set of ``almost'' metric midpoints
between $f(x)$ and $f(y)$.

We claim that $f(U) \subset V_\epsilon$.  Indeed, fix $u\in U$,
then
$$\eqalign{\Vert f(x)-f(u)\Vert &\leq K(\delta )\Vert x-u\Vert \cr
&\leq
{\frac 12} (1+\epsilon) K\Vert x-y\Vert \cr & \leq {\frac 12} (1+
\epsilon)(1-\epsilon)^{-1}\Vert f(x)-f(y)\Vert \cr &\leq (1+4\epsilon)
\Vert f(x) -f(y)\Vert/2 \cr} $$
(because $\epsilon < {\frac 12}$).

On the other hand, if we assume for contradiction that
$$\Vert f(x)-f(u)\Vert \leq (1-4\epsilon)\Vert f(x)-f(y)\Vert /2$$
we obtain
$$\eqalign{(1-\epsilon) K\Vert x-y\Vert &\leq \Vert f(x)-f(y)
\Vert \cr
&\leq \Vert f(x)-f(u)\Vert +\Vert f(u)-f(y)\Vert \cr
&\leq (1-4\epsilon)\Vert f(x)-f(y)\Vert /2 + \Vert f(u)-f(y)\Vert \cr
&\leq K(\delta ) \bigl( (1-4\epsilon)/2 + {\frac 12}\bigr) \Vert
x-y\Vert \cr
&\leq (1+\epsilon) K(1-2\epsilon)\Vert x-y\Vert \ .\cr}$$

So $(1-\epsilon) \leq (1+\epsilon)(1-2\epsilon)$ which is impossible.
Similarly
$$(1-4\epsilon)\Vert f(x)-f(y)\Vert /2 \leq \Vert f(y)-f(u)\Vert \leq
(1+4\epsilon)\Vert f(x)-f(y)\Vert /2 \ .$$
We now specialize to $X=L_1$ and $Y=\ell_1$.

The set $U$ is a ``large'' subset of $L_1$ :\ It contains
an infinite sequence $(x_j)$ so that $\Vert x_j -x_k\Vert =
\Vert x-y\Vert $ for all $j\not= k$.  Indeed, by translation
and change of measure we can assume $x\equiv 0$ and $y\equiv a$ is
constant, and then take $x_j = a(1+r_j)/2$ where  $(r_j)$ are the
Rademacher functions.

On the other hand, $V_\epsilon$ is a ``small'' subset of 
$\ell_1$.  A simple computation shows that $V_\epsilon \subset
C+B_{4\epsilon \Vert f(x)-f(y)\Vert}$, where $C$ is the
compact set of all sequences in $\ell_1$, all of whose coordinates
lie between those of $x$ and $y$.

By the compactness of $C$, there must be $j\not= k$ so that
$\Vert f(x_j)-f(x_k)\Vert \leq 10 \epsilon \Vert f(x)-f(y)\Vert$.

As $\Vert x_j-x_k\Vert =\Vert x-y\Vert  \geq 2\delta > 1$, we 
obtain, using the definition of $L$, that
$$\eqalign{\Vert x-y\Vert  &= \Vert x_j-x_k\Vert \leq L\Vert
f(x_j)-f(x_k) \Vert \cr
&\leq 10\epsilon L \Vert f(x)-f(y)\Vert \leq 10\epsilon L (1+
\epsilon) K\Vert x-y\Vert \ .\cr}$$

So we must have $10\epsilon L(1+\epsilon)K \geq  1$, 
which does not hold if $\epsilon$ is small enough.

\head \S 6.  Non-isomorphic, uniformly-homeomorphic
Banach spaces.\endhead 

In this section we present the recent examples of Ribe \cite{Ri3}
and Aharoni-Lindenstrauss \cite{AL2} of uniformly homeomorphic
Banach spaces which are not isomorphic.

\proclaim{Theorem 6.1}  Let $1\leq p, q, p_n < \infty$ be such
that $p_n\to p$.  Then $(\sum \oplus \ell_{p_n})_q$ is uniformly
homeomorphic to $\ell_p \oplus (\sum \oplus \ell_{p_n})_q$.\endproclaim 

Note that when $p\not= q,p_n \ ,\ \ell_p$ hence also $\ell_p
\oplus (\sum \oplus \ell_{p_n})_q$ does not even Lipschitz
embed into $(\sum \oplus \ell_{p_n})_q$.

As will be evident from the proof the method is quite flexible
and one can prove similar results for other families of spaces,
{\it e.g.},\ $(\sum \oplus L_{p_n})_q$ is uniformly homeomorphic
to $L_p \oplus (\sum \oplus L_{p_n})_q$.

The Theorem was proved for $p=1$ by M.~Ribe \cite{Ri3}.  
I.~Aharoni and J.~Lindenstrauss \cite{AL2} improved Ribe's
techniques so as to yield the Theorem for general $1\leq p<
\infty$, so, in particular, by taking $q,p >1$ we obtain
uniformly convex examples.  Also, by taking $p=1$ and $q$
and $p_n$ all strictly greater than one, we obtain a 
reflexive space uniformly homeomorphic to a non-reflexive
one.

The proof below is a technical simplification of the one
in \cite{AL2}, obtained by carefully making all maps
homogeneous.  This avoids the need to deal simultaneously with the
original norm and with other expressions, uniformly equivalent to
it.

We start with some notation.  In what follows $p$ and $q$ are fixed 
as in the Theorem.  For fixed $1\leq r,s <\infty$ the Mazur
map (see \cite{Maz}) $M_{r,s}\ :\ \ell_r \to \ell_s$ is
defined by
$$\bigl( M_{r,s} (x) \bigr)_n = \Vert x\Vert_r^{1- (r/s)}\ 
x_n^{r/s} \ \hbox{ sign } x_n \ .$$
(Note that we added the normalizing factor $\Vert x\Vert_r^{1-(r/s)}$
so as to have a homogeneous map .)
We shall use the following properties of $M_{r,s}$:
$$\Vert M_{r,s}(x)\Vert_s = \Vert x\Vert_r \ .$$

For each $r,s,t$  $$M_{r,t}=M_{s,t} \circ M_{r,s} \ ,$$ and in 
particular $$M_{r,s}^{-1} = M_{s,r}\ .$$

For each $r,s$ and each $K$,\ $M_{r,s}$ is a uniform 
homeomorphism of the $K$-ball in $\ell_r$ onto the $K$-ball in
$\ell_s$.  Moreover, the  family $\{ M_{r,s}\ :\ 1\leq r,s <
\infty \}$ where each $M_{r,s}$ is restricted to the ball of
radius $\exp (\vert r-s\vert^{-1})$ in $\ell_r$ is equi-uniformly-continuous.
(For a proof of this last fact, see Lemma 1 in \cite{Ri3}.)
\medskip

For each $r$, identify $\ell_r$ with $\ell_r \oplus_r \ell_r$, and
define the map $I_r\ :\ \ell_p \oplus_q \ell_r \to \ell_r = \ell_r
\oplus_r \ell_r$ by
$$I_r (x,y) = {{\Vert (x,y)\Vert}\over {\Vert (M_{p,r}(x),y)\Vert_r}}
\ (M_{p,r}(x),y) \ .$$

Note that even when $r=p$, \ $I_p$ is not the identity but the
renormalization of the norm in $\ell_p \oplus_q \ell_p$ to that of
$\ell_p \oplus_p \ell_p$.  Note also that in this case $I_p$ is a
Lipschitz homeomorphism.

We shall need the following properties of $I_r$:
\smallskip
\itemitem{} $\Vert I_r (a)\Vert =\Vert a\Vert \ .$
\smallskip
\itemitem{} For each fixed $r$ and $K,\ I_r$ is a uniform homeomorphism of
the $K$-ball in $\ell_p \oplus_q \ell_r$ onto the $K$-ball of 
$\ell_r$.
\smallskip
\itemitem{} Moreover, the family $\{ I_r\ :\ 1\leq r<\infty \}$, where each
$I_r$ is restricted to the ball of radius $\exp (\vert p-r\vert^
{-1} )$ in $\ell_p \oplus_q \ell_r$ is equi-uniformly continuous.
\smallskip

To prove the Theorem one can obviously pass to a subsequence of
the $p_n$'s.  We shall thus assume that $\exp (\vert p_n-p\vert^{-1})
> 2^{n+1}$ for each $n=0,1,2,\ldots $\ .

The Theorem follows immediately from the following:

\proclaim{Proposition 6.2}  There is a family of uniform homeomorphisms
$\{ F_t\ :\ 0\leq t<\infty \}$ so that 
\smallskip
\item{(i)}  For $n=1,2,\ldots$ and
$2^{n-1} \leq t\leq 2^n$ $(0\leq t\leq 1$ for $n=0)$, $F_t$ maps the
ball of radius $t$ of $\ell_p \oplus_q \ell_{p_n} \oplus_q \ell_{p_{n+1}}$
onto the ball of radius $t$ of $\ell_{p_n} \oplus_q \ell_{p_{n+1}}$.
\item{(ii)}  $\Vert F_t (a)\Vert =\Vert a\Vert$
\item{(iii)} $F_{2^{n-1}} (x,y,z)=(I_{p_n} (x,y),z))$ for $n=1,2,\ldots $
\item{(iv)} $F_{2^n} (x,y,z)=(y, I_{p_{n+1}} (x,z))$ for $n=0,1,\ldots$
\item{(v)} The family $\{ F_t(a), F_t^{-1}(a)\}$ is equi-uniformly
continuous in both $a$ and $t$, \ {\it i.e.}, there is a function
$\omega (\epsilon) \downarrow 0$ so that
$$\Vert F_t(a)-F_s(b)\Vert ,\  \Vert F_t^{-1}(a)-F_s^{-1}(b)\Vert \leq
\omega (\Vert a-b\Vert +\vert t-s\vert )$$
whenever there is an $n$ s.t. $2^{n-1} \leq t,s\leq 2^n$ and
$a,b$ are such that all expressions make sense.
\smallskip \endproclaim 

(Note that for $t=2^k$, \ $F_t$ is defined twice --- and in different 
ways.  But this should cause no confusion, and in the application
the two definitions will give rise to the same mapping.)

To deduce the Theorem from the Proposition define the
homeomorphism $f :\ \ell_p \oplus_q (\sum \oplus \ell_{p_n} )_q \to
(\sum \oplus \ell_{p_n})_q$ as follows:

If $a=(x,y_0,y_1,\ldots )\in \ell_p \oplus_q (\sum \oplus \ell_{p_n})_q$
is such that $2^{n-1} \leq \Vert a\Vert \leq 2^n$ ($0\leq \Vert
a\Vert \leq 1$ for $n=0$), define
$$f(a)=(y_0, y_1,\ldots , y_{n-1},\ F_{\Vert a\Vert} (x, y_n, y_{n+1}),
\ y_{n+2},\ldots )$$
(with $F_{\Vert a\Vert}(x, y_n, y_{n+1})$ occupying the $n^{th}$
and $(n+1)^{st}$ coordinates).

Given $b =(z_0,z_1,\ldots ) \in (\sum \oplus \ell_{p_n})_q $, 
fixing $n$ so that $2^{n-1} \leq \Vert b\Vert \leq 2^n$ ($n=0$ if
$\Vert b\Vert \leq 1$) determines what are the coordinates to be
changed to obtain a $\epsilon \ \ell_p \oplus_q (\sum \oplus \ell_{p_n})_q$ 
s.t. $f(a)=b$.  Indeed, $a=(x, y_0, \ldots)$ is
given by $y_j =z_j$ if $j\not= n, n+1$, and $(x, y_n, y_{n+1}) =
F_{\Vert b\Vert}^{-1} (z_n, z_{n+1})$.  As all the maps are 
equi-uniformly-continuous, and for $\Vert a\Vert = 2^m$ the two
different formulas obtained by writing $2^{m-1} \leq \Vert a\Vert
\leq 2^m$ or $2^m \leq \Vert a\Vert \leq 2^{m+1}$ agree, $f$ is
indeed a uniform homeomorphism.

\demo{Proof of Proposition}  The proof has several steps.
The first is a formal identification, reducing to a construction in
spaces isomorphic to $\ell_p$.  The second is a simplifying change of
variable. In the third we construct homeomorphisms that do not 
preserve the norm, and in the last step we make the necessary 
normalization.
\medskip
\noindent STEP I.  The maps 
$$g_n\ :\ \ell_p \oplus_q \ell_{p_n} \oplus_q \ell_{p_{n+1}} \longrightarrow
\ell_p \oplus_q \ell_p \oplus_q \ell_p$$
and
$$h_n\ :\ \ell_{p_n} \oplus_q \ell_{p_{n+1}} \longrightarrow \ell_p \oplus_q
\ell_p$$
given by
$$\eqalign{g_n(x,y,z)&= (x, M_{p_n,p}(y), M_{p_{n+1},p}(z))\cr
h_n(x,y)&= (M_{p_n,p}(x), M_{p_{n+1},p}(y))\cr}$$
preserve the norm and are equi-continuous uniform homeomorphisms
when restricted to balls of radius $2^{n+1}$ in the respective spaces.
As $F_t, 2^{n-1} \leq t\leq 2^n$ act on these balls only we shall take
$F_t=h_n^{-1} \circ H_t \circ g_n$ where the maps $$H_t : \ell_p 
\oplus_q \ell_p \oplus_q \ell_p \longrightarrow \ell_p \oplus_q \ell_p$$ will be
constructed to satisfy (i-v) with $p$ replacing both $p_n$ and $p_{n+1}$.
\medskip

\noindent STEP II.  Fixing $n\geq 1$, (only notational changes are
needed for $n=0$), we make a change of variable in $t$, so that it
will vary in the interval $[0,1]$ rather than $2^{n-1} \leq t\leq 2^n$.
Thus we shall construct homeomorphism 
$$G_t\ :\ \ell_p \oplus_q \ell_p \oplus_q \ell_p \longrightarrow \ell_p \oplus_q
\ell_p \quad \hbox{ for } \quad 0\leq t\leq 1$$
satisfying
\smallskip
\item{(a)}  $\Vert G_t (a)\Vert = \Vert a\Vert$, and $G_t$ is homogenous 
\item{(b)}  $G_0 (x,y,z) = (I_p (x,y),z)$
\item{(c)}  $G_1(x,y,z)=(y, I_p(x,z))$
\item{(d)}  There is an absolute constant $K$ so that
$$\Vert G_t(a)-G_s(b)\Vert,\ \Vert G_t^{-1} (a) -G_s^{-1}(a)\Vert
\leq K(\vert t-s\vert +\Vert a-b\Vert )$$
provided $\Vert a\Vert ,\ \Vert b\Vert \leq 1$.
\smallskip

\noindent Once this is done define for $2^{n-1}\leq t\leq 2^n \ \ 
H_t(a)=G_{t2^{1-n}-1}(a)$, and then using the homogeneity of
the $G$'s, we have, whenever $\Vert a\Vert, \ \Vert b\Vert \leq 2^n$
that 
$$\eqalign{\Vert H_t(a)-H_s(b)\Vert &= 2^n\Vert G_{t2^{1-n}-1}
(2^{-n}a)-G_{s2^{1-n}-1}(2^{-n}b)\Vert \cr
&\leq 2^n K (\vert t-s \vert 2^{1-n} +2^{-n} \Vert a-b\Vert )\cr
&\leq 2K (\vert t-s\vert +\Vert a-b\Vert ) \cr}$$
{\it i.e.}, $H_t(a)$ is equi Lipschitzian in both $a$ and $t$, and
similarly for $H_t^{-1}(a)$.
\medskip

\noindent STEP III.  Consider the operators $S_0, S_1\ :\ 
\ell_p \oplus_q \ell_p \oplus_q \ell_p \longrightarrow \ell_p \oplus_q \ell_p$
given by
$$\eqalign{S_0 (x,y,z) &= ( \langle x,y \rangle , z)\cr
S_1 (x,y,z) &= (y, \langle x,z \rangle )\cr}$$
where, in the definition of $S_0$, we identify the first copy of 
$\ell_p$ in $\ell_p \oplus_q \ell_p$ with $\ell_p \oplus_q \ell_p$
and write $\langle x,y \rangle$ for the general point in $\ell_p$
represented this way.  Similarly, in the definition of $S_1$, it
is the second copy of $\ell_p$ in $\ell_p \oplus_q \ell_p$ which
is represented in this way.

The two spaces  $\ell_p \oplus_q \ell_p \oplus_q \ell_p$ and $\ell_p
\oplus_q \ell_p$ are both isomorphic to $\ell_p$, and the general
linear group of $\ell_p$ is contractible.  As both $S_0$ and $S_1$
are isomorphisms, there is a continuous map $t\to S_t$, where for
each $0\leq t\leq 1,\ S_t$ is an invertible operator from
$\ell_p \oplus_q \ell_p \oplus_q \ell_p$ onto $\ell_p \oplus_q
\ell_p$, and a constant $K$, so that 
$S_0$  and  $S_1$  are the given operators and
$$\Vert S_t\Vert,\ \Vert S_t^{-1}\Vert \leq K$$
for all $0 \leq t\leq 1$.

Moreover, using the continuity of $t\to S_t$ and the compactness
of $[0,1]$, we can approximate, and then replace the $S_t$'s
with a map $t\to S_t$ which is Lipschitz in $t$ (in fact, even
piecewise linear).  So we can also assume
$$\Vert S_t -S_s\Vert,\ \Vert S_t^{-1} -S_s^{-1}\Vert \leq
K\vert t-s\vert \ .$$
\medskip

\noindent STEP IV.  The operators $S_t$ almost do the job as
the required $G_t$.  They are, however, only isomorphisms
and we need to preserve the norm, {\it i.e.}, to renormalize
them.  But note that we need special renormalizations for
$t=0$ and $t=1$.  Indeed, $S_0(x,y,z)=(\langle x,y\rangle ,z)$
and we need to have $G_0(x,y,z)=(I_p(x,y),z)$,\ {\it i.e.},
for $t=0$, we need to renormalize the first coordinate of $\ell_p
\oplus_q \ell_p$. 
Similarly, for $t=1$, we need to renormalize the last coordinate
only.  Thus we shall need to do the renormalizations differently
for different values of $t$, and it will be more convenient to work
with the interval $-1\leq t\leq 2$ rather than $0\leq t\leq 1$,
({\it i.e.}, we define $G_t$ for $-1\leq t\leq 2$, and
require in (b) that $G_{-1}(x,y,z)=(I_p(x,y),z)$ and in (c)
that $G_2(x,y,z)=(y, I_p(x,z))$.)

For $0\leq t\leq 1$, we define $G_t(a)=\Vert a\Vert \ {\frac {S_t(a)}
{\Vert S_t(a)\Vert}}$.  Then $G_t^{-1}(b)=\Vert b\Vert {\frac {S_t
^{-1}(b)}{\Vert S_t^{-1}(b)\Vert}}$ and (a) and (d) are satisfied.

Writing $a=(x,y,z)$ we define for $-1\leq t\leq 0$
$$\eqalign{G_t(a)&=\Vert a\Vert {{-t(I_p(x,y),z)+(1+t)S_0(a)}\over
{\Vert -t(I_p(x,y),z)+(1+t)S_0(a)\Vert}} \cr
&= \Vert a\Vert {{(\psi_t(x,y) \langle x,y\rangle ,z)}\over 
{\Vert (\psi_t (x,y) \langle x,y\rangle ,z)\Vert}} \cr}$$
where
$$\psi_t(x,y) = -t {{(\Vert x\Vert^q +\Vert y\Vert^q)^{1/q}}\over
{(\Vert x\Vert^p +\Vert y\Vert^p)^{1/p}}} + (1+t) \ .$$
Then $G_{-1}(x,y,z)=(I_p(x,y),z)$ (because $\Vert (I_p(x,y),z)
\Vert =\Vert a\Vert$), and this formula for $G_0(a)$ agrees with the
previous one, {\it i.e.},
$$G_0(a) =\Vert a\Vert {{S_0(a)}\over {\Vert S_0(a)\Vert}} \ .$$

If $b=(u,v) \in \ell_p \oplus_q \ell_p$ then
$$G_t^{-1} (b)=\Vert b\Vert {{\bigl( {1\over{\psi_t(u)}} u,v\bigr) }\over
{\Vert \bigl( {1\over{\psi_t(u)}} u,v\bigr) \Vert }}$$
(where $\psi_t(u)$ is defined by the above formula upon writing
$u=\langle x,y\rangle$).

It is a straightforward computation to check that $G_t$ 
satisfy also (d) for $-1\leq t\leq 0$.

A similar formula works for $1\leq t\leq 2$:
$$G_t(a) = \Vert a\Vert {{(2-t)S_1(a)+(t-1)(y,I_p (x,z))}\over
{\Vert (2-t) S_1(a)+(t-1)(y, I_p(x,z))\Vert}} \ .$$
\enddemo

\enddocument